\documentclass{amsart}

\usepackage{amsmath}
\usepackage{amssymb}
\usepackage{color}

\usepackage{pdflscape} 
\usepackage{longtable} 

\numberwithin{equation}{section}

\newcommand{\real}{\mathbb{R}}
\newcommand{\nat}{\mathbb{N}}

\newcommand{\li}{\operatorname{li}}
\newcommand{\logpow}[1]{\mathrm{log}^{#1}}

\begin{document}

\title{The combinatorial algorithm for computing $\pi(x)$}

\author[D. B. Staple]{Douglas B. Staple}
\address{Department of Mathematics and Statistics, Dalhousie University}
\email{dstaple@dal.ca}

%
%
%
%
%
%
%
%
%
\subjclass[2010]{Primary 11N05, Secondary 11Y16, 11-04}

\date{\today}

\begin{abstract}
This paper describes recent advances in the combinatorial method for computing $\pi(x)$, the number of primes $\leq x$. In particular, the memory usage has been reduced by a factor of $\log x$, and modifications for shared- and distributed-memory parallelism have been incorporated. The resulting method computes $\pi(x)$ with complexity $O(x^{2/3}\logpow{-2}x)$ in time and $O(x^{1/3}\logpow{2}x)$ in space. The algorithm has been implemented and used to compute $\pi(10^n)$ for $1 \leq n \leq 26$ and $\pi(2^m)$ for $1\leq m \leq 86$. The mathematics presented here is consistent with and builds on that of previous authors.
\end{abstract}

\maketitle

\section{Introduction}

Algorithms used in exact calculations of $\pi(x)$ can be divided into roughly three categories. The simplest algorithms are based on identifying and counting each prime $p\leq x$, typically using some modification of the sieve of Eratosthenes. A na{\"i}ve implementation of the sieve of Eratosthenes uses $O(x\log\log x)$ arithmetic operations and $O(x)$ bits of memory \footnote{When discussing space complexity, one must distinguish between bits of storage and storage locations, each of which grows as the problem size increases. It is commonplace to state that an algorithm has complexity $O(M)$ in space if it requires $\gamma M$ storage locations for some constant $\gamma$, each capable of storing a number with $\log_2 M$ bits \cite{lagariasMCOMP1985,delegliseMCOMP1996,gourdonUnpublished2001b,oliveiraesilvaDETUA2006}. We use this convention here for consistency with other authors. Similarly, in a model of time complexity, one must specify which operations are considered to be performed in constant time. In this paper, we count bitwise operations, addition, subtraction, multiplication, division, modulus, decisions (branches), and memory read and write operations of a single machine word.}. Modern variants based on bucket sieving reduce the memory usage to roughly $\pi(\sqrt{x})$ storage locations each of width $\log_2 \pi(\sqrt{x})$ bits, while leaving the time complexity unchanged \cite{oliveiraesilvaMCOMP2014}. Given the prime number theorem, algorithms that enumerate the primes $p \leq x$ are limited to time complexity $\Omega(x/\log x)$.

The first published algorithm capable of computing $\pi(x)$ substantially faster than the sieve of Eratosthenes was a combinatorial algorithm due to E. Meissel \cite{meisselMA1870}. Given that Meissel's method involved decisions based on human judgement, it is not clear what time complexity to attribute to it; despite this fact, authors usually estimate the time complexity of Meissel's original method as $\Omega(x^{1-\epsilon})$ for any $\epsilon > 0$ \cite{lagariasMCOMP1985}. Meissel used his method in hand calculations of $\pi(10^8)$ and $\pi(10^9)$ in the late 1800s \cite{meisselMA1871,meisselMA1885}; the method was substantially improved by multiple groups of authors, and used in record computations of $\pi(10^n)$ for $10 \leq n \leq 23$ between 1956 and 2007 \cite{lehmerIJM1959,mapesMCOMP1963,bohmanBIT1972,lagariasMCOMP1985,delegliseMCOMP1996,gourdonUnpublished2001b,oliveiraesilvaDETUA2006}. Meissel's method and its descendants are collectively known as ``the'' combinatorial algorithm for computing $\pi(x)$.

Analytic algorithms for computing $\pi(x)$ based on the Riemann zeta function were first presented by Lagarias and Odlyzko in the 1980s \cite{lagariasLNM1984,lagariasJA1987,galwayPhD2004}. Despite the attractive complexity of $O(x^{1/2+\epsilon})$ in time and $O(x^{1/4+\epsilon})$ in space, for any $\epsilon > 0$, the implied constants were large, and no-one succeeded in developing a practical implementation of these methods until nearly 30 years later. The first record computation using an analytic method was $\pi(10^{24})$, under the assumption of the Riemann hypothesis, by Franke, Kleinjung, B{\"u}the, Jost in 2010 \cite{frankeMCOMP2015}. This was followed by a 2012 computation of the same value by Platt without assuming the Riemann hypothesis \cite{plattMCOMP2015}. B{\"u}the \emph{et al.} subsequently modified their algorithm to eliminate the assumption of the Riemann hypothesis, and presented the first computation of $\pi(10^{25})$ \cite{frankeMCOMP2015}.

The problem of determining the number of primes up to some limit is directly tied to the history of the primes themselves, which dates to antiquity and is beyond the scope of this paper. Thus, the paragraphs above are only a sketch of the history of the problem; we direct the interested reader to additional historical references \cite{meisselMA1870,reportBAAS1884,lehmerBook1941,brauerAMM1946,peetreHM1995}.

In the current paper, we describe recent advances to the combinatorial algorithm for computing $\pi(x)$. Firstly, we show how the memory usage of the algorithm can be reduced by a factor of $\log x$. We note that this is not only a reduction in the memory complexity, but a substantial reduction in the actual memory usage for relevant values of $x$. Indeed, before the final step in the memory-complexity reduction was achieved, the author had already reduced the memory usage sufficiently to compute $\pi(10^{26})$, so the original announcement of $\pi(10^{26})$ claimed only a constant-factor reduction in the memory usage. In addition to the reduction in memory usage, we describe mechanisms by which the algorithm can be parallelized. Multiple methods due to the author and others are presented for shared-memory parallelism. We also describe a previously unpublished algorithm for distributed-memory parallelism, loosely based on the idea presented in \cite{gourdonUnpublished2001b}. The algorithms described here were implemented and used to compute $\pi(10^n)$ for $1 \leq n \leq 26$ and $\pi(2^m)$ for $1\leq m \leq 86$.

\section{Reducing space complexity}
\label{section_space}





Three data structures dominate the memory usage in the combinatorial algorithm \cite{lehmerIJM1959,lagariasMCOMP1985,delegliseMCOMP1996,gourdonUnpublished2001b,oliveiraesilvaDETUA2006}: a table of $\pi(y)$ for $y \leq y_{\max}$, a table of the smallest prime factor $p_{\min}(y)$, also for $y \leq y_{\max}$, and a set of $2^L$ sieve counters, where a typical choice for $L\in \nat$ is $L = \left\lfloor \log_2 y_{\max} \right\rfloor$ \cite{delegliseMCOMP1996}. Each of these three data structures limits the space complexity of the algorithm to $O(y_{\max})$. The choice $y_{\max} = \alpha x^{1/3}$ with $\alpha = \beta \logpow{3}x$ for some $\beta \in \real$ is used in the most recent versions of the algorithm \cite{delegliseMCOMP1996,oliveiraesilvaDETUA2006} to achieve the time complexity $O(x^{2/3}\logpow{-2}x)$, simultaneously setting the space complexity at $O(x^{1/3}\logpow{3}x)$. The next largest data structure is a table of primes $p_b$ for $b \leq \pi(y_{\max})$, which has size $O(y_{\max}/\log y_{\max}) = O(x^{1/3} \logpow{2}x)$. Thus, to decrease the memory usage of the algorithm by a factor of $\log x$, we must either reduce each of the limiting data structures by a factor of $\log y_{\max}$ or more, or else eliminate them entirely.

We note that not all expositions of the algorithm are limited by all three of the above data structures. For example, Oliveira e Silva was aware that significantly smaller sieve counters can be used than implied by $L = \left\lfloor \log_2 y_{\max} \right\rfloor$, although he does advocate storing $\pi(y)$ and $p_{\min}(y)$ for $y \leq y_{\max}$ \cite{oliveiraesilvaDETUA2006}. This is to be contrasted with Del{\'e}glise and Rivat, who use $y_{\max}$ sieve counters, and store $\pi(y)$ for $y \leq y_{\max}$, but manage to eliminate $p_{\min}(y)$ from their final formulae \cite{delegliseMCOMP1996}.

\subsection{Retrieving $\pi(y)$ for $y \leq y_{\max}$ in $O(1)$ time using $O(y_{\max}/\log y_{\max})$ space}
\label{section_piy}

The values $\pi(y)$ are used in many places in the algorithm \cite{delegliseMCOMP1996,oliveiraesilvaDETUA2006}. The authors of past studies advocate the use of a table of values for this purpose, which requires $O(y_{\max})$ storage locations. The implied constant is $1$ in the simplest implementation, where a single storage location is used to store a single value of $\pi(y)$. This constant can be reduced somewhat using a wheel, for example only storing $\pi(y)$ for those $y$ coprime to the first $c$ primes, for some $c \in \nat$. However, a wheel cannot be used to reduce the space complexity of the algorithm, as the table used to store the wheel itself grows rapidly, namely with the primorial of $c$. From a practical point of view, even with a wheel the table $\pi(y)$ becomes prohibitively large, and had to be eliminated to permit the computation of $\pi(10^{26})$.

Given the prime number theorem, it turns out that it is possible to retrieve $\pi(y)$ for any $y \leq y_{\max}$ in constant expected time, using only $O(\pi(y_{\max})) = O(x^{1/3} \logpow{2}x)$ precomputed values. The trick is to only store $\pi(\widetilde{y})$ for values $\widetilde{y}$ that are multiples of $\lfloor \log_2 y_{\max} \rfloor$. We also make use of a table of all the primes $p_b$ for $b \leq \pi(y_{\max})$: such a table also requires $\pi(y_{\max})$ storage locations, and is anyway required elsewhere in the combinatorial method \cite{delegliseMCOMP1996,oliveiraesilvaDETUA2006}. The method for determining $\pi(y)$ for a specific value of $y$ is then as follows: firstly, we look up the value $\pi(\widetilde{y})$ at the closest value $\widetilde{y} \leq y$. We then iterate through the array of primes $p_b$, starting at $b = \pi(\widetilde{y}) + 1$, checking whether $p_b > y$ at each value of $b$. If $p_b > y$, we return $\pi(y) = b-1$; if $p_b \leq y$, then we move on to $b+1$, repeating the process.

The surprising thing is the rapid speed with which this algorithm converges: from the prime number theorem, we expect on average one prime in the range $\left( \widetilde{y}, y \right]$, because $y - \widetilde{y} < \lfloor \log_2 y_{\max} \rfloor$. Thus, the most likely situation is that $\pi(\widetilde{y}) = \pi(y)$, i.e., the initial guess for $\pi(y)$ is in fact the correct value, and the algorithm terminates after a single iteration. In practice, in the combinatorial algorithm we retrieve $\pi(y)$ for many values of $y$, such that the average performance is indeed the relevant quantity. Even in the worst case, it is impossible for this algorithm to require more than $\lfloor \log_2 y_{\max} \rfloor$ iterations, which is $O(\log x)$, because this would contradict the assumption that $\widetilde{y}$ was the closest value $\widetilde{y} \leq y$ in the table $\pi(\widetilde{y})$.

\subsection{Iterating over the squarefree $y \leq y_{\max}$ coprime to the first $b$ primes}
\label{section_pmin}

Demanding fast access to $p_{\min}(y)$ for any $y \leq y_{\max}$ is equivalent to factoring any such value of $y$ on demand. $p_{\min}(y)$ is accessed sufficiently often that trivial algorithms such as trial factoring are too slow for this purpose.

The author of \cite{oliveiraesilvaDETUA2006} actually advocated storing the values $p_{\min}(y)\mu(y)$ for $y \leq y_{\max}$, where $\mu(y)$ is the M{\"o}bius function, rather than storing $p_{\min}(y)$ in isolation. However, whether $p_{\min}(y)$ and $\mu(y)$ are stored separately or as a product is immaterial for the current analysis. The values $p_{\min}(y)$ require an array of $y$ storage locations, each of width at least $\log_2 y_{\max}$; the space required to store $\mu(y)$ is negligible by comparison.

As was the case with the array $\pi(y)$, a wheel can be used to compress the array $p_{\min}(y)$. Indeed, the calculation for $\pi(10^{26})$ was performed using a wheel to compress $p_{\min}(y)$ and $\mu(y)$, see Appendix \ref{appendix_benchmarks}. However, even with a wheel the array $p_{\min}(y)$ eventually becomes prohibitively large, and precluded the computation of $\pi(10^{27})$. Luckily, it turns out that the data structure $p_{\min}(y)$ can be completely eliminated, and $\mu(y)$ along with it. In order to do this, we investigate the purpose of storing $p_{\min}(y)$ \cite{oliveiraesilvaDETUA2006}. In fact, the only situation where this array is used is to iterate over all squarefree values $y \leq y_{\max}$ having $p_{\min}(y) > p_b$ for different values of $b \leq \pi(y_{\max})$. The author of \cite{oliveiraesilvaDETUA2006} does this by iterating over all $y \leq y_{\max}$, and explicitly checking the condition $p_{\min}(y) > p_b$ for the given value of $b$. We also note that $\mu(y)$ is used for exactly the same values of $y$. 

Thus, in order to eliminate the array $p_{\min}(y)$, we require an iteration scheme over the squarefree numbers $y \leq y_{\max}$ coprime to the first $b$ primes. Although somewhat cumbersome, it is straightforward to construct such an iteration scheme using a variable number of nested loops. Firstly, we loop over the primes $p_{b_1} \leq y_{\max}$, where $b_1$ is the only loop variable, and assumes the values $[b+1, \pi(y_{\max})]$. We then loop over the biprime numbers $p_{b_1}p_{b_2} \leq y_{\max}$, where $b_1$ ranges from $b+1$ until the product $p_{b_1}p_{b_1+1}$ exceeds $y_{\max}$, and $b_2$ ranges from $b_1 + 1$ until the product $p_{b_1}p_{b_2}$ exceeds $y_{\max}$. We subsequently loop over all numbers $y$ that are the product of three distinct primes $p_{b_1}$, $p_{b_2}$, and $p_{b_3}$, each having $b < b_1 < b_2 < b_3$ and $p_{b_1} p_{b_2} p_{b_3} \leq y_{\max}$, using similar break conditions as above.  This process is repeated until the largest possible number of factors for $y$ has been exceeded, which occurs when $p_{b+1}p_{b+2}p_{b+3} \ldots p_{b+n} > y_{\max}$, where $n$ is the number of nested loops. For example, $p_1p_2\ldots p_{16} = 2\cdot 3\ldots 53 > 2^{64}$, so if $y_{\max}$ is 64 bits or smaller, then $n < 16$. Furthermore, each value of $y = p_{b_1} p_{b_2} \ldots p_{b_n}$ is squarefree by construction, so $\mu(y) = (-1)^n$ for each $y$.

\subsection{Reducing the size of the sieve counters}
\label{section_cnt}

Reducing the size of the sieve counters is easy in comparison to $\pi(y)$ and $p_{\min}(y)$. Firstly, we note that one can simply reduce the number of counters, without negative effects on the runtime \cite{oliveiraesilvaDETUA2006}. By definition, the width of the sieving intervals in the combinatorial algorithm for computing $\pi(x)$ is equal to the number of sieve counters, which we have denoted $2^L$. Given that the upper limit of the sieve is $x / y_{\max}$, there are a total of $x / (2^L y_{\max})$ intervals. Supposing that the overhead per sieving interval is proportional to the number of sieving primes, $\pi(y_{\max})$, the total overhead associated with subdividing the sieving intervals is proportional to $x / (2^L \log x)$ by the prime number theorem. If the overall time complexity is to be kept at $O(x^{2/3}\logpow{-2}x)$, then this implies $2^L > \gamma x^{1/3} \log x$, for some constant $\gamma \in \real$. Choosing this minimal value of $L$ results in sieve counters a factor of $\log x$ smaller than needed to achieve our target space complexity of $O(x^{1/3} \logpow{2}x)$. This is consistent with numerical experiments, where we find that the optimal value of $2^L$ to minimize the runtime is substantially smaller than $y_{\max}$.

Despite the $\logpow{2}x$ reduction above, there is still an incentive to further reduce the size of the sieve counters. Although it is not necessary, it is helpful in a shared-memory architecture to allocate separate sieve counters for each parallel thread. This permits parallelization at the level of sieving blocks, which is sufficiently coarse as to carry relatively little overhead, yet sufficiently fine that load balancing is relatively easy. If such an approach is taken, then the memory usage of the counters is multiplied by a factor of the number of threads $N$, which limits $N$ to $\log x$ or smaller if the memory usage is to be kept at $O(x^{1/3} \logpow{2}x)$.

Two additional approaches for reducing the size of the sieve counters are apparent to the author. Firstly, it should be possible to substantially reduce the amount of overhead per interval using a variant of the bucket sieve algorithm developed by Oliveira e Silva \cite{oliveiraesilvaMCOMP2014}. The basic idea of bucket sieving is to not sieve every interval by every sieving prime, but rather to allocate each sieving prime to a ``bucket'' that indicates the next interval in which a multiple of the prime appears. Buckets are then sequentially processed, one bucket per interval, with each sieving prime encountered being moved to a later bucket. In this fashion, the only primes that are encountered in each sieving interval are the ones for which multiples actually appear in that interval. This permits significantly smaller sieving intervals to be used, effectively eliminating the width of the sieving interval as a contributor to memory usage. Such an approach may even permit the entire sieve table to be stored in the processor's data cache, providing greatly enhanced performance as compared to main memory \cite{oliveiraesilvaMCOMP2014}.

The other potential approach for further reducing the memory usage of the sieve counters involves more efficiently packing the values. The sieve counters suggested by Oliveira e Silva, and used by the present author, have a fractal-like structure \cite{oliveiraesilvaDETUA2006}. For a complete description of the workings and necessity of the sieve counters, we direct the reader to \cite{oliveiraesilvaDETUA2006}. What matters for us is that the counters are each initialized with a number $2^\ell$, for some $0\leq \ell \leq L$, and then decremented from that initial value. This implies that the largest counters need to be stored using integer data types with at least $L$ bits. Thus, if a common binary representation is used for each of the $2^L$ sieve counters, then the total storage requirement is $L2^L$ bits. With the sieve counters indexed using a single variable as in \cite{oliveiraesilvaDETUA2006}, one can probably not avoid using a common binary representation for each of the $2^L$ counters. We note, however, that it is possible to pack the values much more efficiently, resulting in an average of $2$ bits per counter, such that the total requirement is $2^{L+1}$ bits.

\section{Modifications for shared-memory parallelism}
\label{section_onNode}

There are several practical approaches for parallelizing the algorithm on a shared-memory architecture. Firstly, there is one important part of the algorithm, namely the ``easy leaves'' \cite{oliveiraesilvaDETUA2006} in the computation of the partial sieve function $\phi(x, a)$, which can be made embarrassingly parallel. Here $\phi(x, a)$ denotes the count of natural numbers $\leq x$ that are coprime to the first $a$ primes.  The so-called easy leaves do not depend on the main sieve, do not need to be interleaved with other parts of the algorithm, and can be computed completely in isolation of one another.

The difficult part of the parallelism is the main sieve, where the partial sieve function $\phi(m, b)$ is made available for each $m \leq {x/y_{\max}}$ and each prime $b \leq \pi(y_{\max})$. The values of $\phi(m, b)$ for smaller values of $m$ and $b$ are needed in order to compute $\phi(m, b)$ for larger $m$ and $b$, which precludes the embarrassingly parallel computation of $\phi(m, b)$. The approach taken by the current author is to exploit the fact that the sieving is already broken into blocks of length $2^L$. Specifically, one sieves each of $N$ subsequent blocks in parallel, working not with $\phi(m, b)$, but with $\phi(m, b) - \phi(m_{\min}, b)$, where $m_{\min}$ is the beginning of the sieving interval under consideration. Each time a value $\phi(m, b)$ needs to be added to a running sum without knowledge of $\phi(m_{\min}, b)$, this discrepancy is recorded in a tally.  Once each thread is done sieving the interval $[m_{\min}, m_{\min} + 2^L)$, the values $\phi(m_{\min} + 2^L, b) - \phi(m_{\min}, b)$ can be used to compute each $\phi(m_{\min}, b)$, starting at the smallest value of $m_{\min}$, and the discrepancies represented by the tallies can be resolved.

An algorithm that relies on the above idea has several drawbacks. Firstly, separate sieve counters are needed for each thread, which multiplies the memory usage of the sieve counters by a factor of $N$. Secondly, the tallies needed to keep track of the discrepancies between $\phi(m, b)$ and $\phi(m, b) - \phi(m_{\min}, b)$ require a similar amount of memory as the sieve counters. Finally, synchronization is required after each thread sieves a single block, which carries unnecessary overhead. Nonetheless, this approach was found to be efficient enough for the purposes of the author.

After completing the bulk of the current project, the author was made aware of the yet-unpublished work of Kim Walisch. Walisch employs an adaptive algorithm for shared-memory parallelism, where blocks are scheduled dynamically depending on the runtime of previous blocks. Such an approach is certainly more efficient than synchronizing each iteration, which is important if a large shared-memory machine is to be used.

Another potentially attractive approach for shared-memory parallelism, in terms of both time and space, would be to combine adaptive scheduling with the distributed-memory parallelism algorithm that will be explained in the next section.  By leveraging a distributed-memory algorithm even on a shared-memory architecture, the dependence between subsequent iterations would be broken, completely eliminating the need for communication between threads.  Any constant arrays, such as the table of primes $p_b$ for $b \leq \pi(y_{\max})$, could still be shared between the threads to save space on a single shared-memory node.

\section{An algorithm permitting distributed-memory parallelism}
\label{sect_betweenNode}

Distributing the computation of $\pi(x)$ between multiple compute nodes was necessary for the author to compute $\pi(10^{26})$.  The principal issue with distributing the computation is that the simplest algorithms described in Section \ref{section_onNode} rely on rapid exchange of information between compute nodes. Although it is in principle possible to efficiently distribute such a calculation, the greatest degree of parallelism can only be achieved if internode communication can be minimized or eliminated.

Fortunately, it is possible to parallelize the combinatorial algorithm for computing $\pi(x)$ in a way that requires no interprocess communication whatsoever, with the exception of summing the contribution to $\pi(x)$ for each job after the fact. This is highly efficient for the machine, but requires use of a supporting algorithm to break the interdependence of the jobs.

The following algorithm for distributed-memory parallelism is loosely based on an unpublished idea of X. Gourdon \cite{gourdonUnpublished2001b}. Specifically, the issue is that the sums in the main part of the combinatorial algorithm depend on the partial sieve function $\phi(m, b)$, which represents the count of numbers up to $m$ that are coprime to the first $b$ primes. Sieving an interval $[m_{\min}, m_{\min} + 2^L)$ only reveals the values $\phi(m, b) - \phi(m_{\min}, b)$. Thus, determining $\phi(m, b)$ requires storing $\phi(m_{\min}, b)$, updating it after sieving each block, and using the updated value while sieving the next block to obtain any values $\phi(m, b)$ of interest. This approach works fine if the sieve is started at $m_{\min}=0$, because the recursive dependence terminates with $\phi(0, b) = 0$. If the sieve is to be started somewhere in the middle because, for example, earlier blocks are being simultaneously sieved on some other computer, then we need a method to independently compute $\phi(m_{\min}, b)$.

What is needed is an algorithm that can compute $\phi(m, b)$ for a given value of $m = m_{\min}$ and every $c \leq b \leq \pi(y_{\max})$. An idea for how to do this was given in \cite{gourdonUnpublished2001b}, namely to repeatedly apply the recurrence
\begin{equation}
\phi(m, b) = \phi(m, b-1) - \phi(m/p_b, b-1). \label{eq_recurrence}
\end{equation}
Here $c$ is the size of the wheel being used in the sieve, so $\phi(m, c)$ is accessible for any $m \in \nat$ in $O(1)$ time \cite{oliveiraesilvaDETUA2006}. Given $\phi(m, c)$, the idea is to compute $\phi(m/p_{c}, c)$ to obtain $\phi(m, c+1)$. This can be done using the same implementation intended for the overall computation of $\pi(x)$, which is able to compute $\phi(x, a)$ for varying values of $x$ and $a$. The process is then repeated, to obtain $\phi(m, c+2)$, $\phi(m, c+3)$ and onwards up to $\phi(m, \pi(y_{\max}))$.

The difficulty with the above idea is the amount of time needed to perform this process; it would not affect the overall computational complexity of computing $\pi(x)$, but a simple interpretation of this idea was too slow to be used for the computation of $\pi(10^{26})$. The general idea, however, is sound, and modifications can be made to substantially decrease the cost.

The approach taken here is a multifaceted one, where varying methods are used to compute $\phi(m/p_{b}, b)$ depending on the values of $b$. Again, $\phi(m, c)$ is available in $O(1)$ time for any $m \in \nat$ using the sieving wheel. The wheel can also be used to compute $\phi(m, c+1)$ in $O(1)$ time via $\phi(m, c)$ and $\phi(m/p_{c}, c)$. The difficult cases occur for $c+2 \leq b \leq \pi(\sqrt{m})$. We first check whether $p_{b-1}^2 \leq m/p_b$. If this is the case, then we directly apply \eqref{eq_recurrence}, using the combinatorial algorithm to compute $\phi(m/p_b, b-1)$. If, on the other hand, $p_{b-1}^2 > m/p_b$, then Legendre's formula applies, such that $\phi(m/p_b, b-1) = \pi(m/p_b) - b + 2$. We next check whether $m/p_b < y_{\max}$. If this is the case, then we can use the method described in Section \ref{section_piy} to retrieve $\pi(m/p_b)$ in $O(1)$ time. If $m/p_b \geq y_{\max}$ then Legendre's formula still applies, but we must compute $\pi(m/p_b)$ by some other method, e.g., using a second application of Legendre's formula or the combinatorial algorithm. For the remaining values $\pi(\sqrt{m}) < b \leq \pi(y_{\max})$, determining $\phi(m, b)$ is trivial given $\phi(m, b-1)$. Specifically, if $m < y_{\max}$ then $\phi(m, b) = \phi(m, b-1) - 1$ for $\pi(\sqrt{m}) + 1 \leq b \leq \pi(m)$ and $\phi(m, b) = 1$ for $\pi(m) + 1 \leq b \leq \pi(y_{\max})$. If $m \geq y_{\max}$, then $\phi(m, b) = \phi(m, b-1) - 1$ for all $\pi(\sqrt{m}) + 1 \leq b \leq \pi(y_{\max})$.

\section{Numerical results}

\begin{table}[tbp]
\caption{Values of $\pi(x)$ for $x=10^n$}
\label{table_powersOfTen}
\begin{tabular}{r|r|r}
$x$      & $\pi(x)$ & $\li(x) - \pi(x)$\\ \hline \rule{0pt}{1em}
$10^{ 1}$ & 4 & 2.166\\
$10^{ 2}$ & 25 & 5.126\\
$10^{ 3}$ & 168 & 9.610\\
$10^{ 4}$ & 1229 & 17.137\\
$10^{ 5}$ & 9592 & 37.809\\
$10^{ 6}$ & 78498 & 129.549\\
$10^{ 7}$ & 664579 & 339.405\\
$10^{ 8}$ & 5761455 & 754.375\\
$10^{ 9}$ & 50847534 & 1700.957\\
$10^{10}$ & 455052511 & 3103.587\\
$10^{11}$ & 4118054813 & 11587.622\\
$10^{12}$ & 37607912018 & 38262.805\\
$10^{13}$ & 346065536839 & 108971.050\\
$10^{14}$ & 3204941750802 & 314889.954\\
$10^{15}$ & 29844570422669 & 1052618.581\\
$10^{16}$ & 279238341033925 & 3214631.793\\
$10^{17}$ & 2623557157654233 & 7956588.778\\
$10^{18}$ & 24739954287740860 & 21949555.022\\
$10^{19}$ & 234057667276344607 & 99877775.223\\
$10^{20}$ & 2220819602560918840 & 222744643.548\\
$10^{21}$ & 21127269486018731928 & 597394254.333\\
$10^{22}$ & 201467286689315906290 & 1932355208.151\\
$10^{23}$ & 1925320391606803968923 & 7250186215.780\\
$10^{24}$ & 18435599767349200867866 & 17146907278.151\\
$10^{25}$ & 176846309399143769411680 & 55160980939.379\\
$10^{26}$ & 1699246750872437141327603 & 155891678120.791\\
\end{tabular}
\end{table}

\begin{table}[tbp]
\caption{Values of $\pi(x)$ for $x=2^m$}
\label{table_powersOfTwo}
\begin{tabular}{r|r||r|r}
$x$      & $\pi(x)$ & $x$      & $\pi(x)$ \\ \hline \rule{0pt}{1em}
$2^{ 1}$ & 1 & $2^{44}$ & 597116381732\\
$2^{ 2}$ & 2 & $2^{45}$ & 1166746786182\\
$2^{ 3}$ & 4 & $2^{46}$ & 2280998753949\\
$2^{ 4}$ & 6 & $2^{47}$ & 4461632979717\\
$2^{ 5}$ & 11 & $2^{48}$ & 8731188863470\\
$2^{ 6}$ & 18 & $2^{49}$ & 17094432576778\\
$2^{ 7}$ & 31 & $2^{50}$ & 33483379603407\\
$2^{ 8}$ & 54 & $2^{51}$ & 65612899915304\\
$2^{ 9}$ & 97 & $2^{52}$ & 128625503610475\\
$2^{10}$ & 172 & $2^{53}$ & 252252704148404\\
$2^{11}$ & 309 & $2^{54}$ & 494890204904784\\
$2^{12}$ & 564 & $2^{55}$ & 971269945245201\\
$2^{13}$ & 1028 & $2^{56}$ & 1906879381028850\\
$2^{14}$ & 1900 & $2^{57}$ & 3745011184713964\\
$2^{15}$ & 3512 & $2^{58}$ & 7357400267843990\\
$2^{16}$ & 6542 & $2^{59}$ & 14458792895301660\\
$2^{17}$ & 12251 & $2^{60}$ & 28423094496953330\\
$2^{18}$ & 23000 & $2^{61}$ & 55890484045084135\\
$2^{19}$ & 43390 & $2^{62}$ & 109932807585469973\\
$2^{20}$ & 82025 & $2^{63}$ & 216289611853439384\\
$2^{21}$ & 155611 & $2^{64}$ & 425656284035217743\\
$2^{22}$ & 295947 & $2^{65}$ & 837903145466607212\\
$2^{23}$ & 564163 & $2^{66}$ & 1649819700464785589\\
$2^{24}$ & 1077871 & $2^{67}$ & 3249254387052557215\\
$2^{25}$ & 2063689 & $2^{68}$ & 6400771597544937806\\
$2^{26}$ & 3957809 & $2^{69}$ & 12611864618760352880\\
$2^{27}$ & 7603553 & $2^{70}$ & 24855455363362685793\\
$2^{28}$ & 14630843 & $2^{71}$ & 48995571600129458363\\
$2^{29}$ & 28192750 & $2^{72}$ & 96601075195075186855\\
$2^{30}$ & 54400028 & $2^{73}$ & 190499823401327905601\\
$2^{31}$ & 105097565 & $2^{74}$ & 375744164937699609596\\
$2^{32}$ & 203280221 & $2^{75}$ & 741263521140740113483\\
$2^{33}$ & 393615806 & $2^{76}$ & 1462626667154509638735\\
$2^{34}$ & 762939111 & $2^{77}$ & 2886507381056867953916\\
$2^{35}$ & 1480206279 & $2^{78}$ & 5697549648954257752872\\
$2^{36}$ & 2874398515 & $2^{79}$ & 11248065615133675809379\\
$2^{37}$ & 5586502348 & $2^{80}$ & 22209558889635384205844\\
$2^{38}$ & 10866266172 & $2^{81}$ & 43860397052947409356492\\
$2^{39}$ & 21151907950 & $2^{82}$ & 86631124695994360074872\\
$2^{40}$ & 41203088796 & $2^{83}$ & 171136408646923240987028\\
$2^{41}$ & 80316571436 & $2^{84}$ & 338124238545210097236684\\
$2^{42}$ & 156661034233 & $2^{85}$ & 668150111666935905701562\\
$2^{43}$ & 305761713237 & $2^{86}$ & 1320486952377516565496055\\
\end{tabular}
\end{table}

The combinatorial algorithm was implemented and used to compute $\pi(10^n)$ for $1 \leq n \leq 26$ and $\pi(2^m)$ for $1\leq m \leq 86$, see Tables \ref{table_powersOfTen} and \ref{table_powersOfTwo}.  The values $\pi(10^n)$ for $1 \leq n \leq 25$ and $\pi(2^m)$ for $1\leq m \leq 80$ were checked and found to be consistent with the work of previous authors \cite{oliveiraesilvaDETUA2006,frankeMCOMP2015}.  We note that the values $\pi(2^m)$ for $m = 77, 78, 79, 80$ were previously computed under the assumption of the Riemann hypothesis \cite{frankeMCOMP2015}, and were apparently never verified unconditionally until this study.  The values $\pi(10^{26})$ and $\pi(2^m)$ for $81 \leq m \leq 86$ were first reported in this study.  These new values were checked in three ways.  First, each new value was computed twice, using separate clusters and differing numerical parameters ($\alpha$, $c$, and $L$).  Second, the values were checked against the logarithmic integral to ensure the results were reasonable.  Third, at the suggestion of Robert Gerbicz, the parities of the new values of $\pi(x)$ were checked and found to be consistent with those computed by Lifchitz using a yet-unpublished algorithm \cite{lifchitzUnpublished2001}.

\section{Summary}

Recent advances in the combinatorial algorithm for computing $\pi(x)$ were presented together with numerical results. Specifically, memory usage has been reduced by a factor of $\log x$, and algorithms for shared- and distributed-memory parallelism have been developed. The resulting algorithm computes $\pi(x)$ using $O(x^{2/3}\logpow{-2}x)$ arithmetic operations and $O(x^{1/3}\logpow{2}x)$ memory locations, each of width proportional to $\log x$. An algorithm for shared memory parallelism appeared previously in the literature \cite{lagariasMCOMP1985}, but not for the most recent versions of the algorithm \cite{delegliseMCOMP1996,oliveiraesilvaDETUA2006}; the basic idea necessary for distributed memory parallelism appeared in an unpublished manuscript \cite{gourdonUnpublished2001b}. The memory reduction presented here appears to be new.  Previously reported values \cite{oliveiraesilvaDETUA2006,frankeMCOMP2015} of $\pi(10^n)$ for $1\leq n \leq 25$ and $\pi(2^m)$ for $1\leq m \leq 80$ were verified; the values $\pi(10^{26})$ and $\pi(2^m)$ for $81 \leq m \leq 86$ were computed and checked in several ways.

We are now in the interesting situation where two different types of algorithms, combinatorial and analytic, are closely matched for practical calculations of $\pi(x)$. If nothing else, this situation gives unprecedented confidence in any numerical results computed consistently using both types of methods, which is currently the case with $\pi(10^n)$ for $1\leq n \leq 25$ and $\pi(2^m)$ for $1\leq m \leq 80$.

\appendix
\section{Implementation details}
\label{appendix_benchmarks}

\begin{table}[tbp]
\caption{Resources usage for computing $\pi(x)$ with $x=10^n$}
\label{table_benchmarks}
\begin{tabular}{c|c|c|c|c}
         & Time                & Memory                & Time               & Memory\\
         & [node s]            & [bytes]               & [node s]           & [bytes] \\ \hline \rule{0pt}{1em}
$x$      & \multicolumn{2}{c|}{Version 2014.\@10.\@19}       & \multicolumn{2}{c}{Version 2015.\@01.\@30} \\ \hline \rule{0pt}{1em}
$10^{15}$ & $1.48\times{}10^{0}$ & $3.81\times{}10^{7\ }$ & $1.18\times{}10^{0}$ & $1.95\times{}10^{7\ }$ \\
$10^{16}$ & $6.07\times{}10^{0}$ & $4.31\times{}10^{7\ }$ & $5.27\times{}10^{0}$ & $2.16\times{}10^{7\ }$ \\
$10^{17}$ & $2.68\times{}10^{1}$ & $5.78\times{}10^{7\ }$ & $2.59\times{}10^{1}$ & $2.71\times{}10^{7\ }$ \\
$10^{18}$ & $1.31\times{}10^{2}$ & $1.69\times{}10^{8\ }$ & $1.08\times{}10^{2}$ & $1.01\times{}10^{8\ }$ \\
$10^{19}$ & $5.83\times{}10^{2}$ & $3.44\times{}10^{8\ }$ & $6.07\times{}10^{2}$ & $1.74\times{}10^{8\ }$ \\
$10^{20}$ & $2.89\times{}10^{3}$ & $1.73\times{}10^{9\ }$ & $2.56\times{}10^{3}$ & $1.27\times{}10^{9\ }$ \\
$10^{21}$ & $1.20\times{}10^{4}$ & $3.22\times{}10^{9\ }$ & $1.04\times{}10^{4}$ & $1.92\times{}10^{9\ }$ \\
$10^{22}$ & $5.06\times{}10^{4}$ & $5.81\times{}10^{9\ }$ & $4.68\times{}10^{4}$ & $2.98\times{}10^{9\ }$ \\
$10^{23}$ & $2.27\times{}10^{5}$ & $1.16\times{}10^{10}$  & $2.17\times{}10^{5}$ & $5.23\times{}10^{9\ }$ \\
$10^{24}$ & $1.07\times{}10^{6}$ & $2.41\times{}10^{10}$  & -- & $1.00\times{}10^{10}$ \\
$10^{25}$ & $5.25\times{}10^{6}$ & $5.16\times{}10^{10}$  & -- & $2.01\times{}10^{10}$ \\
$10^{26}$ & $2.98\times{}10^{7}$ & $1.12\times{}10^{11}$  & -- & $4.16\times{}10^{10}$ \\
\end{tabular}
\end{table}

The description in \cite{oliveiraesilvaDETUA2006} was used as a starting point for the implementation, with the enhancements of Sections \ref{section_space}--\ref{sect_betweenNode} gradually incorporated.  The implementation was written in the C99 programming language, with significant effort devoted to ensuring the correctness of the program.  Fast unit tests were run on a development machine for every committed version of the code, with more extensive unit tests frequently run on the target cluster.  All code was demanded to compile without warning using the GCC 4.9.1 compiler with the default warning level, and to pass static analysis with the Clang Static Analyzer.  Precisions of finite-width data types were artificially reduced to intentionally break the program and identify failure modes.  Unit tests were written covering wide ranges of parameter values, including edge-cases chosen specifically with the intention of breaking the program.  In general, all code was written and checked as strictly as the author was capable at the time of writing.

In Table \ref{table_benchmarks} we show resources usage for computing $\pi(10^n)$ using two different versions of the author's implementation of the combinatorial algorithm.  The first version of the software, 2014.\@10.\@19, was missing the advancement presented in Section \ref{section_pmin}: this is the version of the software used in the original computations of $\pi(10^{26})$ and $\pi(2^m)$ for $81 \leq m \leq 86$.  In this table, time is measured in ``node seconds'', i.e., it is the sum of the actual time spent on all compute nodes for that calculation.  Similarly, memory usage is memory per node.  Here a ``compute node'' was an IBM iDataplex dx360 M4, having a total of 16 CPU cores (2 $\times$ Intel Xeon E5-2670 eight-core 2.60 GHz CPUs) with either 64 or 128 GB RAM (8 GB PC3-12800 ECC RDIMM modules) depending on the requirements of the calculation.  Thus, $2.98\times{}10^{7}$ node s for computing $\pi(10^{26})$ corresponds to roughly $15.1$ CPU core-years.

\section*{Acknowledgements}

The author thanks Karl Dilcher for support, and for suggestions regarding the underlying algorithm, these calculations, and this paper. Calculations were performed on the Guillimin, Briar{\'e}e, and Colosse clusters from McGill University, Universit{\'e} de Montr{\'e}al, and Laval Universit{\'e}, managed by Calcul Qu{\'e}bec and Compute Canada. The operation of these supercomputers is funded by the Canada Foundation for Innovation (CFI), NanoQu{\'e}bec, RMGA, and the Fonds de recherche du Qu{\'e}bec - Nature et technologies (FRQ-NT).

\ 
\bibliographystyle{amsplain}
\bibliography{pix}

\end{document}